\documentstyle[11pt]{article}

 \newcounter{abceqn}

 \newcounter{abcfig}

\newcommand{\be}{\beta}
\newcommand{\na}{\nabla}
\newcommand{\om}{\omega}
\newcommand{\Ga}{\Gamma}

\newcommand{\pa}{\partial}

\newcommand{\e}{\epsilon}

\newcommand{\k}{\kappa}

\newcommand{\dl}{\delta}
\newcommand{\Dl}{\Delta}
\newcommand{\th}{\theta}
\newcommand{\ra}{\rightarrow}

\renewcommand{\theequation}{\thesection.\arabic{equation}}
\newcommand{\eqnsection}[1]{
	\section{#1}
	\setcounter{equation}{0}
	\renewcommand{\theequation}{\thesection.\arabic{equation}}
	\setcounter{figure}{0}
	\renewcommand{\thefigure}{\thesection.\arabic{figure}}
	\setcounter{remark}{0}
	\renewcommand{\theremark}{\thesection.\arabic{remark}}
	\setcounter{theorem}{0}
	\renewcommand{\thetheorem}{\thesection.\arabic{theorem}}
	\setcounter{lemma}{0}
	\renewcommand{\thelemma}{\thesection.\arabic{lemma}}
}

\title{\bf Global Regularity for the Viscous Boussinesq Equations}

\author{ \\ \\ \\ \\ Yanguang\ \ \ \ Li \\  \\  \\ Department of Mathematics 
\\ \\ University of Missouri \\ \\ Columbia, MO 65211 \\ \\ 
Email: cli@math.missouri.edu}

\date{\today}

\renewcommand{\theequation}{\thesection.\arabic{equation}}

\begin{document}
\bibliographystyle{plain}
\maketitle
\begin{abstract}    
The global regularity for the viscous Boussinesq equations is
proved. \\

AMS MSC (2000): 35, 76.
\end{abstract}

\newtheorem{lemma}{Lemma}
\newtheorem{theorem}{Theorem}
\newtheorem{corollary}{Corollary}
\newtheorem{remark}{Remark}
\newtheorem{definition}{Definition}
\newtheorem{proposition}{Proposition}
\newtheorem{assumption}{Assumption}

\newpage
\tableofcontents

%\newpage

%\thispagestyle{empty}

\newpage
\eqnsection{Introduction}

In the numerical studies \cite{PS92} and \cite{ES94}, different 
observations have been reported on the question of finite time 
singularities of solutions to the inviscid Boussinesq equations:
\begin{eqnarray*}
& & \th_t +u \cdot \na \th =0, \\
& & u_t + u \cdot \na u +\na p =\left ( \begin{array}{c} 0 \\ \th
    \end{array} \right ), \\
& & \na \cdot u =0.
\end{eqnarray*}
Here $\th$ is the temperature, $u=(u_1,u_2)$ is the velocity, $p$ is 
the pressure. In \cite{PS92}, Pumir and Siggia observed that the cap 
of a symmetric rising bubble collapses in a finite time. In contrast,
E and Shu in \cite{ES94} reported that the motion of the bubble cap 
is a very unlikely candidate for finite time singularity formation.

In this paper, we prove the global regularity for the viscous 
Boussinesq equations:
\begin{eqnarray}
& & \th_t +u \cdot \na \th = \k \Dl \th , \nonumber \\
& & u_t + u \cdot \na u +\na p =\left ( \begin{array}{c} 0 \\ \th
    \end{array} \right ) +\nu \Dl u, \label{veq}\\
& & \na \cdot u =0. \nonumber
\end{eqnarray}
Here, $\k$ is the heat diffusion constant, $\nu$ is the viscosity
constant.

The approach of the proof follows from that of \cite{Tem77} \cite{CKN82} 
\cite{CF93}. Nevertheless, we follow closely \cite{CF93}. An possible
alternative approach is the Galerkin approximation approach as in 
chapter $9$ of \cite{CF88}.

\eqnsection{Formulation of the Theorem}

Consider Eq.(\ref{veq}) under the initial conditions:
\begin{equation}
u(0)=b,\ \ \th(0)= \be ; 
\label{veqc}
\end{equation}
where $b$ is $C^\infty$ smooth, divergence free, and compactly 
supported; $\be$ is also $C^\infty$ smooth, compactly supported.
Then the question of global regularity is as follows: Does a smooth 
solution to the Cauchy problem (\ref{veq};\ref{veqc}) exists for 
all time ? For any fixed $T$, we call a solution of (\ref{veq};\ref{veqc})
that is in 
\[
L^{\infty}([0,T];L^2(R^2)) \cap L^2([0,T];H_1(R^2)),
\]
a weak solution; a solutions of (\ref{veq};\ref{veqc}) that is in 
\[
L^{\infty}([0,T];H_1(R^2)) \cap L^2([0,T];H_2(R^2)),
\]
a strong solution. It is known that
strong solutions are actually smooth \cite{CF88}. We will prove the following:
\begin{theorem}
The solution to the Cauchy problem (\ref{veq};\ref{veqc}) 
for $\k >0$, $\nu >0$ is strong 
and hence smooth on the time interval [$0,T$].
\label{thm}
\end{theorem}
and the following:
\begin{corollary}
The solution to the Cauchy problem (\ref{veq};\ref{veqc}) for 
$\k >0$, $\nu =0$ satisfies:
\[
u \in L^{\infty}([0,T];H_1(R^2)),\ \ 
\th \in L^{\infty}([0,T];H_1(R^2)) \cap L^2([0,T];H_2(R^2)).
\]
\end{corollary}

\eqnsection{Proof of the Theorem}

Consider the mollified equation:
\begin{eqnarray}
& & \th_t +u^{(\dl)} \cdot \na \th = \k \Dl \th , \nonumber \\
& & u_t + u^{(\dl)} \cdot \na u +\na p =\left ( \begin{array}{c} 0 \\ \th
    \end{array} \right ) +\nu \Dl u, \label{mveq}\\
& & \na \cdot u = 0. \nonumber
\end{eqnarray}
In which,
\[
u^{(\dl)} =\phi_\dl \ast u,\ \ \na \cdot u^{(\dl)} =0,\ \ \phi_\dl =
\dl^{-2} \phi(x/\dl,y/\dl);
\]
where $\phi$ is smooth, compactly supported, and has integral equal to 1;
\[
\bigg ( \phi_\dl \ast u\bigg )(t,x,y) = \int_{R^2} \phi_\dl(x-x_1,y-y_1)
u(t,x_1,y_1)dx_1dy_1.
\]
It is known that solutions to the Cauchy problem (\ref{mveq};\ref{veqc})
are smooth on any time interval [$0,T$] \cite{CF88}.
The proof consists of $a$ $posteriori$ estimates on the mollified equations. 
The goal is to find uniform bounds in the limit $\dl \ra 0$ and pass to 
the limit to get conclusion for the original Cauchy problem 
(\ref{veq};\ref{veqc}). 

From the first equation in (\ref{mveq}), we get
\begin{equation}
\int \th^2(T) + 2\k \int^{T}_0 \int |\na \th |^2 \leq 
\int \be^2. 
\label{pf1}
\end{equation}
Here $\int$ stands for $\int_{R^2} dxdy$, $\int_0^T$ stands for 
$\int_0^T dt$, and $\th(T)=\th(T,x,y)$.
From the second equation in (\ref{mveq}), we get
\begin{equation}
{d \over dt}\int |u|^2 +2 \nu \int |\na u|^2 = 2 \int \th u_2.  
\label{pf2}
\end{equation}
Apply Schwarz inequality to the right hand side of (\ref{pf2}), and use 
inequality (\ref{pf1}), we have
\begin{equation}
\int |u(T)|^2 + 2 \nu \int^T_0 \int |\na u|^2
\leq \Ga_1,
\label{pf3}
\end{equation}
where $\Ga_1$ is a constant depending only on $b$, $\be$ and $T$.
By inequalities (\ref{pf1};\ref{pf3}), we see that solutions to 
the Cauchy problem (\ref{mveq};\ref{veqc}) are uniformly bounded 
for $\dl \ra 0$ in 
\[
L^{\infty}([0,T];L^2(R^2)) \cap L^2([0,T];H_1(R^2)).
\]
The difficult part is to prove the solutions are uniformly bounded 
in 
\[
L^{\infty}([0,T];H_1(R^2)) \cap L^2([0,T];H_2(R^2)).
\]
Differentiating (\ref{mveq}), we get
\begin{eqnarray}
\left [\pa / \pa t + u^{(\dl)} \cdot \na - \k \Dl \right ]\ \th_x &=&
\{\th, u^{(\dl)}_2\}, \label{deq1} \\
\left [\pa / \pa t + u^{(\dl)} \cdot \na - \k \Dl \right ]\ \th_y &=&
-\{\th, u^{(\dl)}_1\}, \label{deq2} \\
\left [\pa / \pa t + u^{(\dl)} \cdot \na - \nu \Dl \right ]\ \om &=& -\th_x +
\{u_1, v^{(\dl)}_1\} + \{u_2, v^{(\dl)}_2\}. \label{deq3} 
\end{eqnarray}
Here, $\om =\pa u_1 / \pa y - \pa u_2 / \pa x $, $\th_x = \pa \th /\pa x$,
$v^{(\dl)}=u - u^{(\dl)}$, and
\[
\{ p,q\}= {\pa p \over \pa x}{\pa q \over \pa y} -
{\pa p \over \pa y}{\pa q \over \pa x}.
\]
From (\ref{deq1};\ref{deq2};\ref{deq3}), we have:
\begin{eqnarray}
{d \over dt}\int \th^2_x + 2\k \int |\na \th_x |^2 &=& 2 \int \th_x
\{ \th,u^{(\dl)}_2 \}, \label{E1} \\
{d \over dt}\int \th^2_y + 2\k \int |\na \th_y |^2 &=& - 2 \int \th_y
\{ \th,u^{(\dl)}_1 \}, \label{E2} \\
{d \over dt}\int \om^2 + 2\nu \int |\na \om |^2 &=& -2 \int \th_x \om +
2 \int \om [ \{ u_1,v^{(\dl)}_1 \} \label{E3} \\
&+& \{ u_2,v^{(\dl)}_2 \} ]. \nonumber
\end{eqnarray}

The key for the proof below is that we can estimate (\ref{E3}) first,
using the inequality (\ref{pf1}). Then, we estimate (\ref{E1};
\ref{E2}), using the estimates to be obtained
for (\ref{E3}). Notice some simple 
facts:
\begin{eqnarray}
& & \int \om^2 = \int |\na u|^2, \label{fact1} \\
& & \int |\na \om |^2 = \int |\na \na u|^2, \label{fact2} \\
& & \int |v^{(\dl)}|^2 \leq A_1 \dl^2 \int \om^2, \label{fact3} \\
& & \int |\na v^{(\dl)}|^2 \leq A_2 \int \om^2, \label{fact4}
\end{eqnarray}
where $A_1, A_2$ are constants. Next, we estimate (\ref{E3}).
By integration by part, we have
\[
\int \om \{ u_1,v^{(\dl)}_1 \} =\int v^{(\dl)}_1 \{ \om, u_1 \},
\ \ \int \om \{ u_2,v^{(\dl)}_2 \} =\int v^{(\dl)}_2 \{ \om, u_2 \}.
\]
Then, by Schwarz inequality, we get
\begin{equation}
2 \bigg | \int \om [ \{ u_1,v^{(\dl)}_1 \} + \{ u_2,v^{(\dl)}_2 \} ]
\bigg | \leq 4 \bigg \{ \e \int |\na \om |^2 +\e^{-1} \int 
|v^{(\dl)}|^2 |\na u|^2 \bigg \}, \label{pf4}
\end{equation}
where $\e$ is a small constant to be chosen later on. Moreover,
\begin{equation}
\int |v^{(\dl)}|^2 |\na u|^2 \leq \bigg [ \int |v^{(\dl)}|^4 \bigg ]^{1/2}
\bigg [ \int |\na u|^4 \bigg ]^{1/2}. \label{pf5}
\end{equation}
By the Gagliardo-Nirenberg inequality:
\[
\bigg [ \int |v^{(\dl)}|^4 \bigg ]^{1/2} \leq C_1 \bigg [ \int |\na 
v^{(\dl)}|^2 \bigg ]^{1/2} \bigg [ \int |v^{(\dl)}|^2 \bigg ]^{1/2},
\]
and identities (\ref{fact3};\ref{fact4}), we get
\begin{equation}
\bigg [ \int |v^{(\dl)}|^4 \bigg ]^{1/2} \leq C_2 \dl \int \om^2.
\label{pf6}
\end{equation}
From now on, $C_j$ $(j=1,2,...)$ are all constants. 
By the same Gagliardo-Nirenberg 
inequality and identities (\ref{fact1};\ref{fact2}), we get
\begin{equation}
\bigg [ \int |\na u|^4 \bigg ]^{1/2} \leq C_3 \bigg [ |\na \om |^2 
\bigg ]^{1/2} \bigg [ \int \om^2 \bigg ]^{1/2}. 
\label{pf7}
\end{equation}
From (\ref{pf5};\ref{pf6};\ref{pf7}) and an algebraic inequality, we get
\begin{equation}
4 \e^{-1} \int |v^{(\dl)}|^2 |\na u|^2 \leq C_4 \bigg \{ \e \int 
|\na \om |^2 +\e^{-3} \dl^2 \bigg [\int \om^2 \bigg ]^3 \bigg \}.
\label{pf8}
\end{equation}
Thus, by (\ref{pf4};\ref{pf8}), we have 
\begin{equation}
2 \bigg | \int \om [ \{ u_1,v^{(\dl)}_1 \} + \{ u_2,v^{(\dl)}_2 \} ]
\bigg | \leq (4+ C_4)\e \int |\na \om |^2 + C_4 \e^{-3}\dl^2
\bigg [\int \om^2 \bigg ]^3.
\label{pf9}
\end{equation}
Finally, we have the estimate for (\ref{E3}):
\begin{eqnarray}
{d \over dt}\int \om^2 + 2\nu \int |\na \om |^2 &\leq& 2 \bigg [ 
\int \th^2_x \bigg ]^{1/2}\bigg [ \int \om^2 \bigg ]^{1/2} 
\label{pf10} \\
&+& (4+C_4)\e \int |\na \om |^2 + C_4 \e^{-3}\dl^2
\bigg [\int \om^2 \bigg ]^3. \nonumber
\end{eqnarray}
Choose 
\[
\e = {\nu \over 4+C_4},
\]
we get 
\begin{equation}
{d \over dt}\int \om^2 + \nu \int |\na \om |^2 \leq 2 \bigg [ 
\int \th^2_x \bigg ]^{1/2}\bigg [ \int \om^2 \bigg ]^{1/2} +
C_5 \nu^{-3}\dl^2 \bigg [\int \om^2 \bigg ]^3.
\label{pf11}
\end{equation}
Notice that the term 
\[
\bigg [ \int \th^2_x \bigg ]^{1/2}
\]
in (\ref{pf11}) is integrable in time by virtue of (\ref{pf1}). 
Solving the above inequality for $\int \om^2$, and choose 
$\dl$ such that
\[
C_5 \nu^{-3}\dl^2 \bigg \{ 2 \bigg [ \int \om^2(0) \bigg ]^{1/2}
+2 \int^T_0 \bigg [\int \th^2_x \bigg ]^{1/2} \bigg \}^4 T \leq 1,
\]
we get 
\begin{equation}
\int \om^2(T) \leq 4 \bigg \{ \bigg [ \int \om^2(0) \bigg ]^{1/2}
+ \bigg ( {1 \over 2} T \k^{-1} \int \th^2(0) \bigg )^{1/2} \bigg \}^2.
\label{pf12}
\end{equation}
By virtue of (\ref{pf12}), we have from (\ref{pf11}) that
\begin{equation}
\int^T_0 \int |\na \om |^2 \leq \Ga_2,
\label{pf13}
\end{equation}
where $\Ga_2$ is a constant depending only on $b$, $\na b$, $\be$,
and $T$. By (\ref{pf12};\ref{pf13}), we see that $u$ is uniformly 
bounded for $\dl \ra 0$ in 
\[
L^{\infty}([0,T];H_1(R^2)) \cap L^2([0,T];H_2(R^2)).
\]

Next we go back to estimate (\ref{E1};\ref{E2}). By Schwarz inequality,
we have
\begin{equation}
{d \over dt} \int |\na \th |^2 + 2 \k \int |\na \na \th |^2 \leq
8 \bigg [ \int |\na \th |^4 \bigg ]^{1/2} \bigg [ \int |\na 
u^{(\dl)} |^2 \bigg ]^{1/2}.
\label{pf14}
\end{equation}
Notice that 
\[
\int |\na u^{(\dl)} |^2 \leq C_6 \int |\na u |^2;
\]
then (\ref{pf12}) implies that 
\begin{equation}
\bigg [ \int |\na u^{(\dl)} |^2 \bigg ]^{1/2} \leq D,
\label{pf15}
\end{equation}
where $D$ is a constant depending on $T$. By the Gagliardo-Nirenberg
inequality:
\[
\bigg [ \int |\na \th |^4 \bigg ]^{1/2} \leq C_7 \bigg [ \int 
|\na \na \th |^2 \bigg ]^{1/2}\bigg [ \int |\na \th |^2 \bigg ]^{1/2},
\]
we get
\[
{d \over dt} \int |\na \th |^2 + 2 \k \int |\na \na \th |^2 \leq
8 D C_8 \bigg [ \e_1 \int |\na \na \th |^2 + \e^{-1}_1 \int |\na \th |^2
\bigg ].
\]
Choose 
\[
8 D C_8 \e_1 = \k,
\]
we get 
\begin{equation}
{d \over dt} \int |\na \th |^2 +  \k \int |\na \na \th |^2 \leq
D_1 \k^{-1} \int |\na \th |^2,
\label{pf16}
\end{equation}
where $D_1$ is a constant depending on $T$. Then, we have
\[
\int |\na \th(T)|^2 \leq D_2 \int |\na \be |^2 , \ \ 
\int^T_0 \int |\na \na \th |^2 \leq \Ga_3;
\]
where $D_2$ and $\Ga_3$ are constants depending on initial data
and $T$ only. Thus $\th$ is also uniformly bounded for $\dl \ra 0$ 
in 
\[
L^{\infty}([0,T];H_1(R^2)) \cap L^2([0,T];H_2(R^2)).
\]

For this sequence of solutions parametrized by $\dl$, there exists 
a subsequence which converges to a solution to the original Cauchy
problem (\ref{veq};\ref{veqc}); moreover, it inherits the bounds, 
hence it is a strong solution. Since strong solutions are actually
smooth and unique, the theorem is proved. $\clubsuit$

Next, we prove the corollary. Notice that the estimate for 
$\int \om^2(T)$ in (\ref{pf12}) is independent of $\nu$; thus,
we can take the limit $\nu \ra 0$ to get solution to the Cauchy 
problem (\ref{veq};\ref{veqc}) for $\k >0$, $\nu =0$
\[
u \in L^{\infty}([0,T];H_1(R^2)).
\]
All the estimates for $\th$ are still true. This completes the proof
of the corollary. $\clubsuit$

\eqnsection{Conclusion}

The global regularity for the viscous Boussinesq equations is proved,
by following the approach in \cite{CF93}.

\newpage
\bibliography{Bou}

\begin{thebibliography}{1}

\bibitem{CKN82}
L.~Caffarelli, R.~Kohn, and L.~Nirenberg.
\newblock Partial regularity of suitable weak solutions of the
  {N}avier-{S}tokes equations.
\newblock {\em Commun. Pure and Appl. Math.}, XXXV:771--831, 1982.

\bibitem{CF93}
P.~Constantin and C.~Fefferman.
\newblock Directions of vorticity and the problem of global regularity for the
  {N}avier-{S}tokes equations.
\newblock {\em Indiana Univ. Math. J.}, 42 (3):775--789, 1993.

\bibitem{CF88}
P.~Constantin and C.~Foias.
\newblock {\em Navier-{S}tokes {E}quations.}
\newblock Chicago Lectures in Mathematics Series, The University of Chicago
  Press, 1988.

\bibitem{ES94}
W.~E and C-W. Shu.
\newblock Small-scale structures in {B}oussinesq convection.
\newblock {\em Phys. Fluids}, 6 (1):49--58, 1994.

\bibitem{PS92}
A.~Pumir and E.~D. Siggia.
\newblock Development of singular solutions to the axisymmetric {E}uler
  equations.
\newblock {\em Phys. Fluids}, 4:1472, 1992.

\bibitem{Tem77}
R.~Temam.
\newblock {\em Navier-{S}tokes {E}quations. {T}heory and {N}umerical
  {A}nalysis}.
\newblock North-Holland, Amsterdam and New York, 1977.

\end{thebibliography}

\end{document}